\documentclass[12pt]{amsart}

\newcommand{\be}{\begin{equation}}
\newcommand{\ee}{\end{equation}}
\newcommand{\bea}{\begin{eqnarray}}
\newcommand{\eea}{\end{eqnarray}}
\newcommand{\beaa}{\begin{eqnarray*}}
\newcommand{\eeaa}{\end{eqnarray*}}

\newcommand{\lm}{\lambda}

\newcommand{\g}{{\bf g}}
\newcommand{\h}{{\bf h}}
\newcommand{\kk}{{\bf k}}
\newcommand{\lra}{\longrightarrow}
\newcommand{\rla}{\longleftarrow}

\advance\textheight by \topskip
\makeatletter

\@addtoreset{equation}{section}

\def\section{\@startsection {section}{1}{\z@}{-3.5ex plus -1ex minus
 -.2ex}{2.3ex plus .2ex}{\large\bf\centering}}
\def\subsection{\@startsection{subsection}{2}{\z@}{-3.25ex plus%
 -1ex minus -.2ex}{1.5ex plus .2ex}{\bf}}
\def\subsubsection{\@startsection{subsubsection}{3}{\z@}{-3.25ex plus%
 -1ex minus -.2ex}{1.5ex plus .2ex}{\sl}}
\makeatother

\begin{document}

 \baselineskip 17pt
\parindent 10pt
\parskip 9pt

\begin{flushright}
math.QA/0205155\\[3mm]
\end{flushright}
\vspace{1cm}
\begin{center}
{\Large {\bf Rational $K$-matrices and representations\\[0.1in] of
twisted Yangians}}\\ \vspace{1cm} {\large N. J. MacKay}
\\
\vspace{3mm} {\em Department of Mathematics,\\ University of York,
\\York YO10 5DD, U.K.\footnote{emails: {\tt nm15@york.ac.uk} }}
\end{center}

\begin{abstract}
\noindent We describe the twisted Yangians $Y(\g,\h)$ which arise
as boundary remnants of Yangians $Y(\g)$ in 1+1D integrable field
theories. We describe and extend our recent construction of the
intertwiners of their representations (the rational boundary $S$-
or `$K$'-matrices) and perform a case-by-case analysis for all
pairs $(\g,\h)$, giving the $\h$-decomposition of
$Y(\g,\h)$-representations where possible.
\end{abstract}

\section{Introduction}

In recent work \cite{macka01} we looked at what happens to Yangian
($Y(\g)$-)invariant field theories on the half line. From the
exact $S$-matrix point of view, we found that the classes of
solutions of the reflection equation, and thereby the admissible
boundary $S$-matrices, were in correspondence with the $\h$ for
which $(\g,\h)$ is a symmetric pair --- that is, for which $G/H$
is a symmetric space. Each class was then naturally parametrized
by (possibly a finite cover of) $G/H$. From the field theory point
of view, looking in particular at the principal ($G$-valued)
chiral field, we found that the classically integrable boundary
conditions were of two classes. In the one most naturally related
to the above, the field was constrained to take values at the
boundary in $H\subset G$ such that $G/H$ was a symmetric space, or
a translation of $H$. In \cite{deliu01} we found that the
surviving remnant of the $Y(\g)$ symmetry predicts precisely the
reflection matrix structure calculated directly in our earlier
paper.

In this paper we extend the results of \cite{deliu01} on the
construction of boundary $S$-matrices using this boundary
symmetry, the `twisted' Yangian $Y(\g,\h)$. In particular, we
investigate all the $(\g,\h)$, case-by-case, and apply our
techniques wherever possible to obtain the intertwiners of
representations of $Y(\g,\h)$ or `$K$-matrices' --- that is, the
boundary $S$-matrices, up to an overall scalar factor.  This
enables us to list certain of the representations of $Y(\g,\h)$
which contain, as $\h$-irreducible components, the fundamental
representations of $\h$ as their highest components.

Our method is the analogue for the boundary/twisted case of the
`tensor product graph' (TPG) for the bulk case. This is a rather
primitive technique, in the sense that it uses basic conditions on
the $Y(\g,\h)$ representations deduced from Wigner-Eckart theorem
considerations to give the spectral decomposition of the
$K$-matrix. It does {\em not} give explicit constructions of the
$Y(\g,\h)$ action, and breaks down in complex cases. This is
precisely what happens in the bulk
--- where, indeed, the only explicit representations of $Y(\g)$
constructed in the general realization below are those of
Drinfeld's original paper \cite{drinf85}, on certain
$\g$-irreducible representations, and on $\g\oplus{\mathbb C}$.

\section{Yangians $Y(\g)$}

Suppose the Lie algebra $\g$ to be generated by $Q_0^a$ with
structure constants $f^a_{\;\;bc}$ and (trivial) coproduct
$\Delta$,
\begin{equation}\label{cr0} \left[ Q_0^a , Q_0^b \right] = i
f^{a}_{\;\;bc} Q_0^c \hspace{0.3in}{\rm and}\hspace{0.3in}
\Delta(Q_0^a) = Q_0^a \otimes 1 + 1 \otimes Q_0^a \,.
\end{equation}
The Yangian \cite{drinf85} $Y(\g)$ is the enveloping algebra
generated by these and $Q_1^a$, where
\begin{equation}\label{cr1}
\left[ Q_0^a , Q_1^b \right] = i  f^{a}_{\;\;bc} Q_1^c
\hspace{0.3in}{\rm and}\hspace{0.3in} \Delta(Q_1^a) = Q_1^a
\otimes 1 + 1 \otimes Q_1^a + {1\over  2}f^{a}_{\;\;bc} Q_0^b
\otimes Q_0^c\,.
\end{equation}
The requirement that $\Delta$ be a homomorphism fixes\footnote{For
$\g\neq sl(2)$. For the general condition see
Drinfeld\cite{drinf85}.}
\begin{equation}\label{YSerre}
 f^{d[ab} [Q_1^{c]},Q_1^d] \; = \; {i\over {12}}
  \,f^{ap i} f^{bq j}f^{cr k}f^{ijk} \, Q_0^{(p} Q_0^q Q_0^{r)} \;,
\end{equation}
where $(\,)$ denotes symmetrization and $[\,]$ anti-symmetrization
on the enclosed indices, and indices have been raised and lowered
freely with the invariant metric $\gamma$.

The Yangian may be considered as a deformation of the polynomial
algebra $\g[z]$: with $Q_1^a=zQ_0^a$, the undeformed algebra would
satisfy (\ref{YSerre}) with the right-hand side zero -- that is,
$z^2$ times the Jacobi identity. In $Y(\g)$, (\ref{YSerre}) acts
as a rigidity condition on the construction of higher $Q^a_n$ from
the $Q_1^a$. There is an (`evaluation') automorphism
$$L_\theta\,:\; Q_0^a \mapsto Q_0^a\,,\qquad Q^a_1 \mapsto  Q_1^a+
\theta { c_A \over 4i\pi} Q^a_0 \,,$$ where $c_A=C_2^\g(\g)$ is
the value of the quadratic Casimir $C_2^\g\equiv \gamma_{ab}Q_0^a
Q_0^b$ in the adjoint representation. (We have chosen this
normalization so that, in integrable quantum field theories with
$Y(\g)$ symmetry, $\theta$ is the particle rapidity.)

Thus any representation $v$ of $Y(\g)$ may be considered as
carrying a parameter $\theta$: the action of $Y(\g)$ on $v^\theta$
is that of $L_\theta(Y(\g))$ on $v^0$. The $i$th fundamental
representation $v_i^\theta$ of $Y(\g)$ is in general reducible as
a $\g$-representation, with one of its irreducible components
(that with the greatest highest weight, where these are partially
ordered using the simple roots) being the $i$th fundamental
representation $V_i$ of $\g$. In the simplest cases (which include
all $i$ for $\g=a_n$ and $c_n$), $v_i^\theta=V_i$ as a
$\g$-representation, and $Q_1^a=\theta { c_A \over 4i\pi} Q^a_0$
upon it.

One way to deduce the $\g$-irreducible components of the other
$v_i$ is to use the fusion procedure: one constructs $v_i^\theta
\otimes v_j^{\theta'}$ using the $\g$-irreducible $v_i,v_j$ of the
last paragraph, and then notes that while $v_i^\theta \otimes
v_j^{\theta'}$ is generally $Y(\g)$-irreducible, for certain
special values of $\theta-\theta'$ it may be $Y(\g)$-reducible
(though not fully reducible) to another fundamental representation
$v_k^0$.

This can be seen using the tensor product graph (TPG)
\cite{macka91,gould02}. One constructs a graph whose nodes are the
$\g$-irreducible components of $v_i^\theta \otimes v_j^{\theta'}$,
with edges between nodes $U$ and $V$ when $Q_1^a$ has non-trivial
action from $U$ to $V$. The edge labels (whose calculation is
detailed in \cite{macka91}) are then the special values at which
reducibility occurs. For example, let $v_1^\theta=V_1$, the vector
representation of $so(N)$. The graph of $v_1^\theta\otimes
v_1^{\theta'}$ is then $$ (2\lm_1)\quad \stackrel{{2i\pi\over
N}}{\lra}\quad (\lm_2)\quad\stackrel{i\pi}{\lra} \quad(0) \;,$$
where we have labeled the representations by their highest
weights, so that $(\lm_i)\equiv V_i$. At $\theta=-\theta'=i\pi/N$,
this becomes reducible: the action of $Q_1^a$ on $V_2$ no longer
yields states in $(2\lm_1)$, and we have constructed $v_2^0$,
which decomposes as a $so(N)$-representation into $(\lm_2)\oplus
(0)$.

These decompositions for general $\g$ and $i$ appeared
incrementally in the literature \cite{ogiev86,chari95}; for a full
enumeration for simply-laced $\g$ see \cite{kleber96}. Many
further results can be deduced using the TPG for non-fundamental
representations: if we remove from any TPG all the edges with a
particular label, the remaining subgraphs each provide
representations of $Y(\g)$. Thus in the above example
$Y(\g)$-representations can be constructed whose
$\g$-decomposition is $(\lm_2)\oplus(0)$, $(2\lm_1)$,
$(2\lm_1)\oplus(\lm_2)$ and $(0)$.

\section{Twisted Yangians $Y(\g,\h)$}

Let $\h\subset\g$ be the subalgebra of $\g$ invariant under an
involutive automorphism $\sigma$. We shall write $\g=\h \oplus
\kk$, so that $\h$ and $\kk$ are the subspaces of $\g$ with
$\sigma$-eigenvalues $+1$ and $-1$ respectively. We shall use
$a,b,c,...$ for general $\g$-indices, $i,j,k,...$ for $\h$-indices
and $p,q,r,...$ for $\kk$-indices.

We define the {\bf twisted Yangian} \cite{deliu01} $Y(\g,\h)$ to
be the subalgebra of $Y(\g)$ generated by \bea\label{Q0q}
&&Q_0^i\\{\rm and} && \widetilde{Q}_1^p \equiv Q_1^p + {1\over 4}
[C, Q_0^p] \label{Q1q}\,, \eea where $C\equiv \gamma_{ij}Q_0^i
Q_0^j$ is the quadratic Casimir operator of $\g$,
 restricted to $\h$.

We can again consider $Y(\g,\h)$ to be a deformation, this time of
the subalgebra of (`twisted') polynomials in $\g[z]$ invariant
under the combined action of $\sigma$ and $z\mapsto -z$. Its
defining feature is that $Y(\g,\h)$ is a co-ideal subalgebra
\cite{deliu01}, $\Delta (Y(\g,\h)\subset Y(\g)\times Y(\g,\h)$.
(In an integrable-model setting, in which $Y(\g,\h)$ is the
symmetry of a model with boundary, this allows the boundary states
to form representations of $Y(\g,\h)$ while the bulk states form
representations of $Y(\g)$.) It specializes to the cases studied
in \cite{molev96}, though the relationship between the two
approaches remains to be fully explored.

The TPG may be generalized to deal with $Y(\g,\h)$. In the same
way that $Y(\g)$-representations were generally $\g$-reducible, so
$Y(\g,\h)$-representations naturally form representations of the
subalgebra $\h\subset Y(\g,\h)$, and these are generally
reducible. The key idea is the branching graph \cite{deliu01}: one
considers $v_i$ as a $\g$-representation and determines how this
reduces further as an $\h$-representation. In the simplest cases,
where $v_i=V_i$ is $\g$-irreducible, the $\h$-irreducible
components of $V_i$ are the nodes of the graph, and the edges
connect those $\h$-irreducible representations (hereafter
`irreps') between which $Q_0^p$ has non-trivial action, while the
labels are constructed from the differences in $C$ between these
components -- for details and many examples see \cite{deliu01}.

The more subtle cases are those where $v_i$ is $\g$-reducible, and
we rederive the branching graph here for this more general case.
First recall that the $Y(\g,\h)$-representations are intertwined
by the $K$- or `reflection' matrix $ K_v(\theta): v^\theta
\rightarrow v^{-\theta}$, where $v^\theta$ is a
$Y(\g)$-representation. (In certain cases, in which
non-self-conjugate $\g$-representations branch to self-conjugate
$\h$-representations, $v^\theta$ may be conjugated by $K$
\cite{macka01}.) Intertwining the $Q_0^i$ (that is, from the
physics point-of-view, their conservation in boundary scattering
processes) requires that $$ K_v(\theta) Q_0^i = Q_0^i K_v(\theta)
$$ (in which by $Q_0^i$ we mean here its representation on $v$)
and thus that $K_v(\theta)$ act as the identity on
$\h$-irreducible components of $v$. So we have $$K_v(\theta) =
\sum_{W \subset V\subset v^\theta} \tau_W(\theta) P_W\,, $$ where
the sum is over $\h$-irreps $W$ into which the $V$ branch, where
$V$ is a $\g$-irreducible component of $v^\theta$; $P_W$ is the
projector onto $W$.

To deduce relations among the $\tau_W$ we intertwine the
$\widetilde{Q}_1^p$. Recall that, {\em within} a $\g$-irreducible
$V\subset v^\theta$, the action of $Q_1^p$ is given by $ Q_1^p =
\theta { c_A \over 4i\pi} Q_0^p$, so that $$\langle W \vert\vert
K_v(\theta) \left( \theta { c_A \over 4i \pi} Q_0^p + {1\over
4}[C, Q_0^p] \right)\vert\vert W' \rangle = \langle W
\vert\vert\left( -\theta { c_A \over 4 i\pi}Q_0^p + {1\over 4}[C,
Q_0^p] \right) K_v(\theta)\vert\vert W' \rangle \,,$$ for
$W,W'\subset V$. Thus when the reduced matrix element $\langle W
\vert\vert Q_0^p\vert\vert W' \rangle\neq 0$ we have \be\label{BG}
{\tau_{W'}(\theta) \over \tau_{W}(\theta)} =
\left[\Delta\right]\,, \qquad{\rm where} \;\left[ A \right] \equiv
{\frac{i\pi A}{c_A}+\theta \over \frac{i\pi A}{c_A} -\theta}\ee
and $\Delta = C(W)-C(W')$. To find the $W,W'$ for which $\langle W
\vert\vert Q_0^p\vert\vert W' \rangle\neq 0$ we recall that $\kk$
forms an irreducible representation $K$ of $\h$ (reducing into two
conjugate representations of $\tilde\h$, with opposite $u(1)$
numbers, where $\h=\tilde\h\times u(1)$). A necessary condition
for (\ref{BG}) to apply is then that $W\subset K\otimes W'$.
Although not automatically sufficient, this is (as in the bulk
case \cite{gould02}) sufficiently constraining in simple cases to
enable us to deduce $K$.

There are also links between $W,W'$ which descend from {\em
different} $\g$-irreps. When $W\subset K\otimes W'$, there will
generally be some (unknown) action of $Q_1^p$ between them, but
since $\langle W \vert\vert Q_0^p\vert\vert W' \rangle = 0$, we
will have $\tau_W=\tau_{W'}$.

We then describe $K_v(\theta)$ by using a graph, in which the
nodes are the equivalence classes of $\h$-irreps quotiented by the
relation $W\sim W' \Leftrightarrow W\subset K\otimes W'$ and $ W
\subset V \neq V' \supset W'$. These classes, the (generally
$\h$-reducible) representations $\widetilde{W}_i$, are linked by
an edge, directed from $\widetilde{W}_i$ to $\widetilde{W}_j$ and
labelled by $\Delta_{ij}$, whenever $W\subset K\otimes W'$ for any
$W\subset \widetilde{W}_i,$ $W'\subset \widetilde{W}_j$. (Note
that the $\Delta_{ij}$ calculated from all such pairs $W,W'$ are
equal.)

To calculate the labels, we first write $$ C = \sum_i c_i
C_2^{\h_i}\,,$$ where $\h=\bigoplus_i \h_i$ is a sum of simple
factors $\h_i$ (and $C_2^{\h_i}$ is the quadratic Casimir of
$\h_i$). We now compute the $c_i$, with $\gamma$ the identity both
on $\g$ and on each $\h_i$, by taking the trace of the adjoint
action\footnote{I should like to thank Tony Sudbery for this
suggestion.} of $C$ on $\g$, yielding $$ c_i= {c_A \over
C_2^{\h_i}(\h_i) + {{\rm dim}\, \kk\; \over {\rm dim}\, \h_i}
C_2^{\h_i}(\kk)}\,.$$ The relative values of the $c_i$ are of
course crucial in the cases where $\h$ is non-simple, but there is
also a highly non-trivial implication of their absolute values:
for the edge of the physical strip for the boundary $S$-matrices
to be at $\theta=i\pi/2$, $$ \sum_i \left( {C_2^{\h_i}(\h_i) \over
C_2^{\h_i}(\kk)} + {{\rm dim}\, \kk\; \over {\rm dim}\,
\h_i}\right)^{-1}={1\over 2}$$ must hold. That it does so, and
does so only for symmetric spaces, is a result of \cite{godd85},
also known as the `symmetric space theorem' \cite{daboul96}.

As an example of a graph, let $v^\theta$ be the $\g$-reducible
example we encountered earlier, $v_2^0=(\lm_2)\oplus(0)$ of
$so(N)$, and let $\h=so(M)\times so(N\!-\!M)$. We denote irreps of
$\h$ by $(\mu,\nu)$, where $\mu$ is an $so(M)$ weight and $\nu$ an
$so(N\!-\!M)$ weight. Then the graph is $$ (0,\lm_2)
\quad\stackrel{N-2M-2}{\lra}\quad (\lm_1,\lm_1)\oplus(0,0)
\quad\stackrel{N-2M+2}{\lra}\quad (\lm_2,0)\,.$$ For $v_3^0 =
(\lm_3)\oplus(\lm_1)$ we find similarly
 $$ (0,\lm_3) \quad\stackrel{N-2M-4}{\lra}\quad
(\lm_1,\lm_2)\oplus(0,\lm_1) \quad\stackrel{N-2M}{\lra}\quad
(\lm_2,\lm_1)\oplus(\lm_1,0)
\quad\stackrel{N-2M+4}{\lra}\quad(\lm_3,0)\,.$$ Thus we see that,
at certain special values of $\theta$, the graph truncates, and
the action of $Y(\g,\h)$ may be consistently restricted: $v_2$ to
$(\lm_2,0)$ or $(0,\lm_2)$, $v_3$ to $(\lm_3,0)$ or $(0,\lm_3)$.

There are various limitations to this method which cause it to
break down in cases more complex than those we shall treat. First,
as commented upon in the introduction, when $W,W'$ branch from the
same $U$, the condition $W\subset K\otimes W'$ is necessary but
not automatically sufficient for (\ref{BG}) to apply. Second, the
method breaks down -- just as does the TPG for the bulk case --
when any $W$ appears with multiplicity greater than one. Third,
when $W,W'$ branch from different $U$, the action of $Q_1^p$ is
generally unknown, and so we do not know precisely for which
$W,W'$ we have $\langle W \vert\vert Q_1^p\vert\vert W' \rangle
\neq 0$.

\section{The rational $K$-matrices}

First, a general result. Recall \cite{drinf85} that
$\g\oplus{\mathbb C}$ always extends to a representation of
$Y(\g)$. This branches to $\h\oplus\kk\oplus{\mathbb C}$ of
(simple) $\h$. As a representation of $Y(\g,\h)$ the branching
graph is $$ \h \lra \kk\oplus {\mathbb C}$$ (with label
$(C(\h)-C(\kk))$), and we see that $\kk\oplus{\mathbb C}$ extends
to a representation of $Y(\g,\h)$. In the cases where $\h$ is not
simple, the graph consists of a central node $\kk\oplus{\mathbb
C}$ and edges between this and each component of $\h$. We can only
construct $\kk\oplus{\mathbb C}$ as a representation of $Y(\g,\h)$
in this way when all the factors of $\h$ are isomorphic.

We now proceed to a case-by-case analysis, detailing all the
$(\g,\h)$ and $v_i$ for which our method yields results. In each
case we give the representation $K$ of $\h$ formed by $\kk$. When
$\h$ contains a $u(1)$ factor, we write $\h=\tilde\h\times u(1)$;
we then give the decomposition of $K$ into irreps of $\tilde{\h}$
(we do not specify the values of the $u(1)$ generator). For
classical $\g$ it is simplest to deal with the $b$- and $d$-
series together as $so(N)$, and so we treat the other classical
groups in the same way, as $su(N)$ and $sp(2n)$. When we move on
to the exceptional cases, we write $a_n=su(n+1)$, $b_n=so(2n+1)$,
$c_n=sp(2n)$ and $d_n=so(2n)$ in the usual way.

We first give the decomposition of the fundamental
$Y(\g)$-representations $v_i$ into fundamental
$\g$-representations $V_i$ (though, for $e_7$ and $e_8$, only for
those $i$ of which we will be able to make use). Next we give the
graphs for $K_{v_i}(\theta)$, where these can be computed. For the
cases in which $\h$ has only one non-trivial simple factor, we
then list the $Y(\g,\h)$-representations which follow, in terms of
$\h$- (or $\tilde\h$-) representations $(\lm)$ given in terms of
their highest weights $\lm$, where the fundamental weights are
$\lm_i$, following the Dynkin diagram conventions in the appendix.
Where the weight of the top component is $\lm_i$, we label these
$w_i$. Throughout, $\lfloor x\rfloor$ denotes the integer part of
$x$.

\subsection{$\g=su(N)$}

$v_i=V_i,\quad$ $i=1,2,\ldots,N\!-\!1$.

\noindent Most of the $su(N)$ cases are contained in
\cite{macka01}, but we include them for completeness.

\subsubsection{$\h=su(M) \times su(N\!-\!M) \times u(1)$ \\
$K=(\lm_1,\lm_{m-1})\oplus(\lm_{n-m-1},\lm_1)$} \label{A} The
graphs for $v_r$, $r=1,\ldots,=\lfloor N/2\rfloor$ (the others
follow by conjugation), are $$
\hspace{-0.1in}(0,\lm_r)\stackrel{N-2M-2(r-1)}{\lra}
(\lm_1,\lm_{r-1})\lra \ldots (\lm_p,\lm_{r-p})
\stackrel{N-2M-2(r-1)+4p}{\lra}\ldots(\lm_{r-1},
\lm_1)\stackrel{N-2M+2(r-1)}{\lra}(\lm_r,0)\,.$$ Note that
$(\lm_M,0)\equiv(0,0)$ and $(0,\lm_{N-M})\equiv(0,0)$, while
$(\lm_r,0)$ vanishes for $r>M$ and $(0,\lm_r)$ vanishes for
$r>N\!-\!M$, causing the graph to truncate. Then
 \beaa w_r & = & (\lm_r,0),\quad
r=1,2,\ldots, M\!-\!1 \\ w'_r & = & (0,\lm_r),\quad r=1,2,\ldots,
N\!-\!M\!-\!1\,.\eeaa Many other (non-fundamental) representations
can be constructed \cite{short02}, and may involve graphs which
are $p$-dimensional (in the sense that they contain nodes linked
by $2p$ edges).

\subsubsection{$\h=so(N)$\\ $K=(2\lm_1)$}

Here there are no non-trivial graphs; we simply have the
$\g\rightarrow\h$ branching rules for $V_r$: \beaa &  \qquad w_r =
V_r\quad{\rm for}\quad r=1,2,\ldots,\lfloor
(N-3)/2\rfloor,\quad{\rm then}&
\\[0.1in] &
\begin{array}{ll}N \quad{\rm even:}\quad & w_{(N-2)/2} =
 (\lm_s+\lm_{s^\prime})\,,\quad
w_{N/2}
 =(2\lm_s)\oplus(2\lm_{s^\prime})\\[0.1in]
N \quad {\rm odd:} & w_{(N-1)/2} = (2\lm_s).\end{array}& \eeaa

\subsubsection{$\h=sp(N),\;$ $N=2n$\\ $K=(\lm_2)$}

The graph for $v_r$ is  $$ (\lm_{r})
\stackrel{N-2(r-2)}{\lra}(\lm_{r-2})\ldots
\stackrel{N-2(r-2p)}{\lra} (\lm_{r-2p})\stackrel{N+4-2(r-2p)}{
\lra}\ldots \left\{\begin{array}{llll} (\lm_{2})\stackrel{N}{\lra}
& (0) & \qquad$r$\;{\rm even}\\(\lm_{3})\stackrel{N-2}{ \lra} &
(\lm_{1}) & \qquad$r$\;{\rm odd}\end{array}\right.\,,$$so that
$$w_r = (\lm_r)\oplus(\lm_{r-2})\oplus\,\ldots\,\oplus
\left\{\begin{array}{ll} (0) & r\;{\rm even}\\(\lm_1) & r\;{\rm
odd} \end{array}\right.$$

\subsection{$\g=so(N)$}

This is the only classical case for which the $v_i$ are generally
reducible as $so(N)$-representations: $v_i=V_i\oplus V_{i-2}
\oplus \,\ldots\, V_{0/1}$ for $r=1,2,\ldots,[(N-3)/2]$,
$v_s=V_s$, $v_{s'}=V_{s'}$.

\subsubsection{$\h=so(M) \times so(N\!-\!M)$\\ $K=(\lm_1,\lm_1)$}

For $r=1,2,\ldots,[(N-3)/2]$, the graphs are
  $$
\hspace{-0.1in}(0,r)\stackrel{N-2M-2(r-1)}{\lra} (1,{r-1})\lra
\ldots (p,{r-p}) \stackrel{N-2M-2(r-1)+4p}{\lra}\ldots({r-1},
1)\stackrel{N-2M+2(r-1)}{\lra}(r,0)$$ (as in the $su(N)$ case),
where (for $p<q$; $p\geq q$ is analogous)
$(p,q)\equiv(\lm_p,\lm_q)\oplus(\lm_{p-1},\lm_{q-1})
\oplus\,\ldots\,(0,\lm_{q-p})$, and so \beaa w_r & = &
(\lm_r,0),\quad r=1,2,\ldots, [(M-3)/2] \\ w'_r & = &
(0,\lm_r),\quad r=1,2,\ldots, [(N-M-3)/2]\,.\eeaa For $v_s$ and
$v_{s'}$ we have $$ v_{s^{(\prime)}}=w_{s^{(\prime)}}=
(\lm_s,\lm_{s^\prime})\oplus(\lm_{s^\prime},\lm_s)$$ if $V_s\neq
V_{s^\prime}$ for either $so(M)$ or $so(N\!-\!M)$;
$v_s=(\lm_s,\lm_s)$ if not.

\subsubsection{$\h=su(n)\times u(1),\;$ $N=2n$\\
$K=(\lm_2)\oplus(\lm_{n-2})$}

We deal first with the spinor representations, distinguishing the
cases $n=2m+1$ (in which $V^*_s=V_{s^\prime}$, where $*$ denotes
complex conjugation) from $n=2m$ (in which $V^*_{s^{(\prime)}} =
V_{s^{(\prime})}$). The graphs for $V_s$ and $V_{s^\prime}$ are
respectively (noting $\lm_n=0$ and $(\lm_{n-r})=(\lm_r)^*$) $$
(\lm_{2m}) \stackrel{2(n+3-4m)}{\lra}
\;\ldots\;\stackrel{2(n-1-4q)}{\lra}
(\lm_{2q})\stackrel{2(n+3-4q)}{\lra}\ldots
\lra(\lm_4)\stackrel{2(n-5)}{\lra}(\lm_2)
\stackrel{2(n-1)}{\lra}(0) $$ and
$$(\lm_1)\stackrel{2(5-n)}{\lra}(\lm_3) \lra
\;\ldots\;\stackrel{2(4q+1-n)}{\lra}(\lm_{2q+1})\stackrel{2(4q+5-n)}
{\lra}\ldots \left\{\begin{array}{ll}
\stackrel{2(n-3)}{\lra}(\lm_{2m-1}) & \quad n=2m
\\\stackrel{2(n-1)}{\lra} (\lm_{2m+1}) &\quad  n=2m\!+\!1
\end{array}\right.\,.$$ Thus we have, for $n=2m+1$, \beaa w_{2p}&
= & \bigoplus_{r=0}^p (\lm_{2r})
\\ w_{2p+1} & = &\bigoplus_{r=0}^p (\lm_{2r+1})\,,\eeaa  and
their conjugates (for $p=0,1,2,\ldots,m$); while, for $n=2m$,
\beaa w_{2p} & = & \bigoplus_{r=0}^p (\lm_{2r}) \oplus
(\lm_{n-2r})\,,\quad p=0,1,\ldots,[(m-1)/2] \\ w_{2p+1} & =
&\bigoplus_{r=0}^p (\lm_{2r+1}) \oplus (\lm_{n-2r-1})\,,\quad
p=0,1,\ldots,[(m-2)/2] \\ w_m & = & (\lm_m) \oplus w_{m-2}
\,,\eeaa which are all self-conjugate.

Turning to the antisymmetric tensor representations, $v_1=V_1$
branches trivially to $(\lm_1)\oplus(\lm_{n-1})$; then for
$v_2=V_2\oplus 1$ the graph is $$ (\lm_1+\lm_{n-1})
\stackrel{2}{\lra} (\lm_2)\oplus(\lm_{n-2})\oplus(0)
\stackrel{2n-2}{\lra} (0)\,.$$ Thereafter, for $i=3,4,\ldots,n-2$,
the graph becomes intractable for the reasons mentioned at the end
of the last section.

\subsection{$\g=sp(2n) $}

$v_i=V_i,\quad$ $i=1,2,\ldots,n$.

\subsubsection{$\h=sp(2m) \times sp(2n\!-\!2m)$ \\ $K=(\lm_1,\lm_1)$}

The graphs are as in section \ref{A}, and the $w_r$ are therefore
\beaa w_r &=& (\lm_r,1),\quad r=1,2,\ldots, m\\ w'_r & = &
(1,\lm_r),\quad r=1,2,\ldots, n-m\,.\eeaa

\subsubsection{$\h=su(n) \times u(1)$ \\
$K=(2\lm_1)\oplus(2\lm_{n-1})$}

The $sp(2n)\rightarrow su(n)$ branching rule is $$ V_r =
\bigoplus_{a=0}^r (\lm_a+\lm_{n-r+a})\,, $$ from which the graph
is$$ (\lm_r) \stackrel{2-2r}{\lra}(\lm_{r-1}+\lm_{n-1})
\lra\ldots\stackrel{4a-2-2r}{\lra}
(\lm_{r-a}+\lm_{n-a})\stackrel{4a+2-2r}{\lra}\ldots
(\lm_1+\lm_{n-r+1}) \stackrel{2r-2}{\lra}(\lm_{n-r}).$$ We thus
have $$ w_r =(\lm_r)\oplus(\lm_{n-r})\,,\qquad
r=1,2,\ldots,[n/2]\,, $$

\subsection{$\g=e_6$}

$v_1=V_1,\, v_6=V_6;\; v_2=V_2\oplus {\mathbb C}$

\subsubsection{$\h=c_4$\\ $K=(\lm_4)$}

\beaa v_1 \;{\rm or}\; v_6:\qquad && (\lm_2)\\[0.05in] v_2: \qquad
&& (\lm_4)\oplus(0)\stackrel{2}{\lra}(2\lm_1)\\[0.05in] v_3 \;{\rm
or}\;v_5: \qquad && (\lm_1+\lm_3)\oplus(\lm_2)
\stackrel{6}{\lra}(2\lm_1) \eeaa so that
$w_2=(\lm_2),\,w_4=(\lm_4)\oplus(0)$.

\subsubsection{$\h=d_5 \times u(1)$\\
$K=(\lm_4)\oplus (\lm_5)$}

\beaa v_1:
\qquad&&(\lm_1)\stackrel{2}{\rla}(\lm_5)\stackrel{10}{\lra}
(0)\\[0.05in] v_2 : \qquad&& (\lm_2)\stackrel{4}{\lra}
(\lm_4)\oplus(\lm_5)\oplus (0) \stackrel{12}{\lra} (0)\eeaa and
$v_6$ as $v_1^*$; the others are intractable.\\ Hence
$w_1=(\lm_1),\,
w_2=(\lm_2),\,w_4=(\lm_4)\oplus(0),\,w_5=(\lm_5)\oplus(0)$, others
unknown.

\subsubsection{$\h=a_5\times a_1$\\
$K=(\lm_3,\lm_1)$}

Here we abbreviate the (spin-$s/2$) $su(2)$ irrep $(s\lm_1)$ to
$(s)$. \beaa v_1:\qquad && (\lm_4,0)\stackrel{2}{\lra}
(\lm_1,1)\\[0.05in] v_2 : \qquad && (\lm_5+\lm_1,0)
\stackrel{0}{\lra} (\lm_3,1)\oplus(0,0) \stackrel{8}{\lra}
(0,2)\,,\eeaa $v_6$ as $v_1^*$, and we can go no further.

\subsubsection{$\h=f_4$\\ $K=(\lm_4)$}

\beaa v_1:\qquad && (\lm_4)\stackrel{12}{\lra} (0)
\\[0.05in] v_2: \qquad&& (\lm_1)\stackrel{6}{\lra}(\lm_4)\oplus(0)
\eeaa and $w_1=(\lm_1)$ and $w_4=(\lm_4)$ in addition to
$\kk\oplus{\mathbb C}=(\lm_4)\oplus(0)$.

\subsection{$\g=e_7$}

$v_1=V_1\oplus{\mathbb C},\; v_7=V_7,\;v_2=V_2\oplus V_7$ are the
only cases we can treat.

\subsubsection{$\h=e_6 \times u(1)$\\ $K=(\lm_1)\oplus(\lm_6)$}
$$ v_1:\qquad (\lm_2)\stackrel{6}{\lra}
(\lm_1)\oplus(\lm_6)\oplus(0)\stackrel{18}{\lra}(0)$$ and
$w_2=(\lm_2)$.

\subsubsection{$\h=d_6 \times a_1$\\ $K=(\lm_5,1)$}
\beaa v_7:\qquad&& (\lm_6,0)\stackrel{4}{\lra}(\lm_1,1)\\[0.05in]
 v_1:\qquad&&
(\lm_2,0)\stackrel{2}{\lra}(\lm_5,1)\oplus(0,0)\stackrel{14}{\lra}(0,2)
\eeaa

\subsubsection{$\h=a_7$\\ $K=(\lm_4)$}
\beaa v_7:\qquad&& (\lm_2)\oplus(\lm_6) \\[0.05in]
 v_1:\qquad &&(\lm_4)\oplus(0)\stackrel{2}{\lra} (\lm_1+\lm_7)\\[0.05in]
 v_2:\qquad &&(2\lm_1)\stackrel{8}{\rla}
 (\lm_1+\lm_5)\oplus(\lm_2)\oplus
 (\lm_3+\lm_7)\oplus(\lm_6)\stackrel{8}{\lra}(2\lm_7)\eeaa
and $w_4=(\lm_4)\oplus(0)\,(=\kk+{\mathbb C})$.

\subsection{$\g=e_8$}

$v_8=V_8\oplus{\mathbb C}$

\subsubsection{$\h=e_7 \times a_1$\\ $K=(\lm_7,1)$}

$$v_8:\qquad (\lm_1,0)\stackrel{6}{\lra}
(\lm_7,1)\oplus(0,0)\stackrel{26}{\lra} (0,0)$$ and
$w_1=(\lm_1,0)$.

\subsubsection{$\h=d_8$\\ $K=(\lm_7)$}
\beaa v_8:\qquad && (\lm_7)\oplus(0)\stackrel{2}{\lra} (\lm_2)
\\[0.05in] v_1:\qquad && (\lm_4)\oplus(\lm_7)\oplus(0)
\stackrel{2}{\lra} (\lm_1+\lm_8)\oplus(\lm_2) \stackrel{14}{\lra}
(2\lm_1) \eeaa and $w_2=(\lm_2)$.

\subsection{$\g=f_4$}

$v_4=V_4,\, v_1=V_1\oplus{\mathbb C}$

\subsubsection{$\h=b_4$ \\ $K=(\lm_4)$}
\beaa v_4:\qquad &&
(\lm_1)\stackrel{1}{\rla}(\lm_4)\stackrel{9}{\lra}(0)
\\[0.05in] v_1: \qquad && (\lm_2)\stackrel{5}{\lra} (\lm_4)\oplus(0) \eeaa so
$w_1=(\lm_1),\,w_2=(\lm_2),\, w_4=(\lm_4)\oplus(0)$.

\subsubsection{$\h=c_3 \times a_1$ \\ $K=(\lm_3,1)$}
\beaa v_4:\qquad && (\lm_2,0)\stackrel{1}{\lra} (\lm_1,1)
\\[0.05in] v_1:\qquad && (2\lm_1,0) \stackrel{1}{\rla}
 (\lm_3,1)\oplus(0,0)\stackrel{5}{\lra}
(0,2)\eeaa

\subsection{$\g=g_2$}
$v_1=V_1,\, v_2=V_2\oplus{\mathbb C}$

\subsubsection{$\h=a_1 \times a_1$ \\ $K=(3,1)$}
\beaa v_1:\qquad && (2,0)\stackrel{-2/3}{\lra}(1,1) \\ v_2:\qquad
&& (2,0) \stackrel{8/3}{\rla}
(3,1)\oplus(0,0)\stackrel{0}{\lra}(0,2)\eeaa

\pagebreak
\section{$K$-matrices and the magic square}

The Freudenthal-Tits magic square (see \cite{barto02} and
references therein) provides a remarkable construction, based on
division algebras, of the exceptional Lie algebras. We recall here
only that it may be written $$
\begin{array}{c|cccc}
m= & \quad 1 & 2 & 4\quad & 8 \\ \hline  &\quad  a_1 & a_2 &
c_3\quad  & f_4
\\ & \quad a_2 & a_2 \times a_2 & a_5 \quad & e_6 \\
 & \quad c_3 & a_5 & d_6\quad  & e_7
\\ &\quad  f_4 & e_6 & e_7\quad  & e_8
\end{array}$$
It is no coincidence that, if we write $g_m^{(i)}$ for the
appropriate entry of the $i$th row, then
$g_m^{(4)}/(g_m^{(3)}\times a_1)$, $g_m^{(3)}/(g_m^{(2)}\times
u(1))$ and $g_m^{(2)}/g_m^{(1)}$ are all symmetric spaces --- this
is fundamental in the construction of the square. The parameter
$m$ is the order of an underlying division algebra (real or
complex numbers, quaternions or octonions, here as a derivation
algebra; the rows are labelled in the same way by triality
algebras). For each of $i=1,2,3$ the dimension of the
corresponding representation $K$ is a linear function of $m$,
respectively $3m+2$, $2(3m+3)$ and $2(6m+8)$.

It has already been noted \cite{westb02} that the $R$-matrices
(the solutions of the bulk Yang-Baxter equation) in these
distinguished representations have a uniform graph  structure for
each $i$, with the graph labels having a simple linear dependence
on $m$. We note here that the same is true of the $K$-matrices we
have constructed: for $g^{(2)}_m/g^{(1)}_m$, there are
$K$-matrices \beaa & \h \stackrel{2m-4}{\lra}\kk\oplus{\mathbb C}
&\\& \kk \stackrel{3m}{\lra} {\mathbb C},& \eeaa \noindent while
for $g^{(3)}_m/(g^{(2)}_m\times u(1))$ we have $$ \h
\stackrel{m-2}{\lra} \kk\oplus{\mathbb C}\stackrel{2m+2}{\lra}{\bf
C}.$$ \noindent Finally for $g^{(4)}_m/(g^{(3)}_m\times a_1)$
there is $$(\h,0) \stackrel{m-2}{\lra}(\kk,1)\oplus{\mathbb C}
\stackrel{3m+2}{\lra}(0,2)$$ and, for all but $m=8$, $$ (U,0)
\stackrel{m}{\lra}(V,1)\,,$$ where $V$ is the vector
representation of $g^{(3)}_m$ and $U$ a representation of which we
have no general characterization. This even extends, just as in
the bulk case, to $g^{(4)}_0=d_4$ and $g^{(4)}_{-2/3}=g_2$.

At this level, of course, the above is merely a nice observation,
but it does suggest that it might be interesting to study Yangians
and twisted Yangians -- indeed, Yang-Baxter and reflection
equation algebras more generally -- from a division algebra point
of view.

\noindent{\bf Acknowledgments}

I should like to thank Tony Sudbery and Bruce Westbury for
discussions, and Ben Short, Max Nazarov and Gustav Delius for
helpful comments.

 \pagebreak
\section{APPENDIX: Dynkin diagrams and conventions}

$su(N)=a_n$, $N=n+1$  :

\begin{picture}(250,20)
\put(20,15){\circle{6}} \put(50,15){\circle{6}}
\put(80,15){\circle{6}} \put(140,15){\circle{6}}
\put(170,15){\circle{6}} \put(200,15){\circle{6}}
\put(23,15){\line(1,0){24}} \put(53,15){\line(1,0){24}}
\put(143,15){\line(1,0){24}} \put(173,15){\line(1,0){24}}
\put(105,15){$\dots$} \put(18,0){$1$} \put(48,0){$2$}
\put(78,0){$3$} \put(126,0){$n\!-\!2$} \put(160,0){$n\!-\!1$}
\put(198,0){$n$}
\end{picture}

\vspace{0.1in} \noindent $so(N)=b_n$, $N=2n+1$  :

\begin{picture}(250,20)
\put(20,15){\circle{6}} \put(50,15){\circle{6}}
\put(80,15){\circle{6}} \put(140,15){\circle{6}}
\put(170,15){\circle{6}} \put(200,15){\circle*{6}}
\put(23,15){\line(1,0){24}} \put(53,15){\line(1,0){24}}
\put(143,15){\line(1,0){24}} \put(172,13){\line(1,0){26}}
\put(172,17){\line(1,0){26}} \put(181,12){$>$}
\put(105,15){$\dots$} \put(18,0){$1$} \put(48,0){$2$}
\put(78,0){$3$} \put(126,0){$n\!-\!2$} \put(160,0){$n\!-\!1$}
\put(198,0){$s$}
\end{picture}

\vspace{0.1in} \noindent $sp(2n)=c_n$  :

\begin{picture}(250,20)
\put(20,15){\circle*{6}} \put(50,15){\circle*{6}}
\put(80,15){\circle*{6}} \put(140,15){\circle*{6}}
\put(170,15){\circle*{6}} \put(200,15){\circle{6}}
\put(23,15){\line(1,0){24}} \put(53,15){\line(1,0){24}}
\put(143,15){\line(1,0){24}} \put(172,13){\line(1,0){26}}
\put(172,17){\line(1,0){26}} \put(181,12){$<$}
\put(105,15){$\dots$} \put(18,0){$1$} \put(48,0){$2$}
\put(78,0){$3$} \put(126,0){$n\!-\!2$} \put(160,0){$n\!-\!1$}
\put(198,0){$n$}
\end{picture}

\vspace{0.1in} \noindent $so(N)=d_n$, $N=2n$  :

\begin{picture}(250,35)(0,-10)
\put(20,15){\circle{6}} \put(50,15){\circle{6}}
\put(80,15){\circle{6}} \put(140,15){\circle{6}}
\put(170,15){\circle{6}} \put(191,36){\circle{6}}
\put(191,-6){\circle{6}} \put(23,15){\line(1,0){24}}
\put(53,15){\line(1,0){24}} \put(143,15){\line(1,0){24}}
\put(172,17){\line(1,1){17}} \put(172,13){\line(1,-1){17}}
\put(105,15){$\dots$} \put(18,0){$1$} \put(48,0){$2$}
\put(78,0){$3$} \put(152,-2){$n\!-\!2$} \put(200,33){$s$}
\put(200,-9){$s'$}
\end{picture}

\vspace{0.1in} \noindent $e_6$  :

\begin{picture}(250,25)
\put(20,15){\circle{6}} \put(50,15){\circle{6}}
\put(80,15){\circle{6}} \put(110,15){\circle{6}}
\put(140,15){\circle{6}} \put(80,45){\circle{6}}
\put(23,15){\line(1,0){24}} \put(53,15){\line(1,0){24}}
\put(83,15){\line(1,0){24}} \put(113,15){\line(1,0){24}}
\put(80,18){\line(0,1){24}} \put(18,0){$1$} \put(48,0){$3$}
\put(78,0){$4$} \put(108,0){$5$} \put(138,0){$6$} \put(86,42){$2$}
\end{picture}

\vspace{0.1in} \noindent $e_7$  :

\begin{picture}(250,30)
\put(20,15){\circle{6}} \put(50,15){\circle{6}}
\put(80,15){\circle{6}} \put(110,15){\circle{6}}
\put(140,15){\circle{6}} \put(170,15){\circle{6}}
\put(80,45){\circle{6}} \put(23,15){\line(1,0){24}}
\put(53,15){\line(1,0){24}} \put(83,15){\line(1,0){24}}
\put(113,15){\line(1,0){24}} \put(143,15){\line(1,0){24}}
\put(80,18){\line(0,1){24}} \put(18,0){$1$} \put(48,0){$3$}
\put(78,0){$4$} \put(108,0){$5$} \put(138,0){$6$} \put(168,0){$7$}
\put(86,42){$2$}
\end{picture}

\vspace{0.1in} \noindent $e_8$  :

\begin{picture}(250,30)
\put(20,15){\circle{6}} \put(50,15){\circle{6}}
\put(80,15){\circle{6}} \put(110,15){\circle{6}}
\put(140,15){\circle{6}} \put(170,15){\circle{6}}
\put(200,15){\circle{6}} \put(80,45){\circle{6}}
\put(23,15){\line(1,0){24}} \put(53,15){\line(1,0){24}}
\put(83,15){\line(1,0){24}} \put(113,15){\line(1,0){24}}
\put(143,15){\line(1,0){24}} \put(173,15){\line(1,0){24}}
\put(80,18){\line(0,1){24}} \put(18,0){$1$} \put(48,0){$3$}
\put(78,0){$4$} \put(108,0){$5$} \put(138,0){$6$} \put(168,0){$7$}
\put(198,0){$8$} \put(86,42){$2$}
\end{picture}

\vspace{0.1in} \noindent $f_4$  :

\begin{picture}(250,10)(-20,0)
\put(20,15){\circle{6}} \put(50,15){\circle{6}}
\put(80,15){\circle*{6}} \put(110,15){\circle*{6}}
\put(23,15){\line(1,0){24}} \put(52,17){\line(1,0){26}}
\put(52,13){\line(1,0){26}} \put(83,15){\line(1,0){24}}
\put(61,12){$>$} \put(18,0){$1$} \put(47,0){$2$} \put(78,0){$3$}
\put(108,0){$4$}
\end{picture}

\noindent $g_2$  :

\begin{picture}(250,10)(-40,0)
\put(20,15){\circle*{6}} \put(50,15){\circle{6}}
\put(22,17){\line(1,0){26}} \put(23,15){\line(1,0){24}}
\put(22,13){\line(1,0){26}} \put(30,12){$<$} \put(18,0){$1$}
\put(48,0){$2$}
\end{picture}


\parskip 8pt

\begin{thebibliography}{99}

\bibitem{macka01} N. MacKay and B. Short, {\em Boundary scattering,
symmetric spaces and the principal chiral model on the half-line},
Comm. Math. Phys. {\bf 233}(2003)313, {\tt hep-th/0104212}

\bibitem{deliu01}
G. Delius, N. MacKay and B. Short, {\em Boundary remnant of
Yangian symmetry and the structure of rational reflection
matrices}, Phys. Lett. {\bf B522}(2001)335, {\tt hep-th/0109115}

\bibitem{drinf85}
V. Drinfeld, {\em Hopf algebras and the quantum Yang-Baxter
equation}, Sov. Math. Dokl. {\bf 32}(1985)254

\bibitem{macka91}
N. J. MacKay, {\em Rational R-matrices in irreducible
representations}, J. Phys. {\bf A24}(1991)4017

\bibitem{gould02}
M. Gould and Y-Z. Zhang, {\em R-matrices and the tensor product
graph method}, {\tt hep-th/0205071} and references therein

\bibitem{godd85}
P. Goddard, W. Nahm and D. Olive, {\em Symmetric spaces,
Sugawara's energy momentum tensor in two dimensions and free
fermions}, Phys. Lett. {\bf B160}(1985)111

\bibitem{daboul96}
C. Daboul, {\em Algebraic proof of the symmetric space theorem},
J. Math. Phys. {\bf 37}(1996)3576, {\tt hep-th/9604108}

\bibitem{ogiev86}
E. Ogievetsky, N. Reshetikhin and P. Wiegmann, {\em The principal
chiral field in two dimensions on classical Lie algebras}, Nucl.
Phys. {\bf B280} (1987) 45 ; E. Ogievetsky and P. Wiegmann, {\em
Factorized $S$-matrix and the Betha ansatz for simple Lie groups},
Phys. Lett. {\bf B168}(1986)360

\bibitem{chari95}
V. Chari and A. Pressley, {\em Fundamental representations of
Yangians and singularities of $R$-matrices}, J. reine angew. Math.
{\bf 417}(1991)87

\bibitem{kleber96}
M. Kleber, {\em Combinatorial structure of finite dimensional
representations of Yangians: the simply-laced case}, Int. Math.
Res. Notices {\bf 4}(1996)187, {\tt q-alg/9711032}

\bibitem{molev96}
A. Molev, M. Nazarov and G. Ol'shanskii, {\em Yangians and
classical Lie algebras}, Russ. Math. Surveys {\bf 51}(1996)205; A.
Molev, {\em Finite-dimensional irreducible representations of
twisted Yangians}, J. Math. Phys. {\bf 39}(1998)5559, {\tt
q-alg/9711022} ; A. Molev and E. Ragoucy, {\em Representations of
reflection algebras}, {\tt math.QA/0107213};
 E. Ragoucy, {\em
Quantum group symmetry of integrable systems with or without
boundary}, {\tt math.QA/0202095}

\bibitem{short02}
B. Short, PhD thesis, unpublished.

\bibitem{barto02}
C. Barton and A. Sudbery, {\em Magic squares and matrix models of
Lie algebras}, {\tt math.RA/0203010}

\bibitem{westb02}
B. Westbury, {\em R-matrices and the magic square}, J. Phys. {\bf
A}36(2003)2857

\end{thebibliography}
\end{document}